\documentclass[12pt]{article}
\usepackage{amssymb}
\usepackage{amsmath,amsfonts,amsthm}
\usepackage{srcltx}
\usepackage{footnote}
\usepackage{tikz-cd}
\begin{document}                                 
	   \begin{center} \bf{IDEAL CONTAINMENT  Vs. POWERS}\\ Pramod K. Sharma\\ e-mail: pksharma1944@yahoo.com    \end{center}
	   
	  \begin{center}    ABSTRACT 
	   
	\end{center}

 Let $R$ be a commutative ring with identity. In this note, we study the property: If $ I \subsetneqq  J$ are ideals in $R$, then $ I^n \subsetneqq J^n$ for all $ n\geq 1$. We define the notion of a big ideal (Definition 1.2). It is noted that the property has close relationship with the notions of reduction of an ideal and Ratliff-Rush ideal [7]. Apart from other results, it is proved that a Noetherian domain satifies the property if and only if every ideal in $R$ is a Ratliff-Rush ideal. We also prove that  ideals having no proper reduction are big ideals, and maximal ideals in regular rings are big.

  \section{INTRODUCTION}
  
	Throughout this note, all rings are commutative with identity $ (\neq 0).$  While working on [10], we needed to know which rings satisfy  (P) : Whenever $ I \subsetneqq J$ are ideals then $ I^n \subsetneqq J^n$ for all $n \geqq 1$. Absence of any information on this question is the reason for this note. We define (Definition 1.2) an ideal $J$ in a ring $R$ to be a big ideal if whenever  $I\subsetneqq J$ then $ I^n \subsetneqq J^n$ for all $n \geqq 1$. Thus a ring $R$ satisfies the property (P) if every ideal in $R$ is a big ideal.
	         
	         In section 2, apart from some general results, we prove that a Noetherian integral domain $R$ satisfies the property (P) if and only if every ideal in $R$ is a Ratliff-Rush ideal (Definition 4.2)., and also prove that if a Noetherian integral domain $R$ satisfies the property (P), then dimension of $R$ is $ \leqq 1$. Further, we show that a Dedekind domain satisfies the property (P).

	         The section 3 deals with the existance of big ideals. We prove that an ideal $J$ in a Noetherian interal domain $R$ is a big ideal if and only if whenever $ I \subsetneqq J$  is an ideal then $ J \subsetneqq I^*$. Further, if an ideal $J$ in a ring $R$ admits no proper reduction, then $J$ is a big ideal. We also prove that any maximal ideal in a  regular ring is a big ideal.

	      \section{Preliminaries on Rings Satisfying (P)}  We shall study here the following property :
	         (P): A ring $R$ satisfies the property if whenever  $I\subsetneqq J$ are ideals in  $R$  then   $ I^n \subsetneqq J^n$ for all $n \geqq 1$.

	         In this connection, we define:\\	
	         
	         {\bfseries Definition} 1.2.  An ideal $J$ in a ring $R$ is called a big ideal if whenever an ideal $ I \subsetneqq J$ , then $ I^n \subsetneqq J^n$ for all $ n \geqq 1$.\\ 
	         
	         {\bfseries Remark} 2.2.  A ring $R$ satisfies the property (P) if every ideal in $R$ is a big ideal.\\
	         
	          First of all, we record some  definitions for convenience of the reader.\\
	         
	        {\bfseries Definition} 3.2. Let $I \subset J$ be ideals in a ring $R$. We shall say thst $I$ is a reduction of $J$ if $ IJ^m = J^{m+1}$ for some $ m \geq 1$.\\
	         
	         {\bfseries Definition} 4.2. Let $I$ be an ideal in a Noetherian ring $R$, then $I^* = \bigcup \{ (I^{m+1} : I^m) : m  \geq 1 \} $ . If $I$ is regular i.e, $I$ contains a non-zero divisor, and $I = I^*$ , then $I$ is called a RATLIFF-RUSH ideal of $R$.\\ 
	         
	           \textbf{Lemma} 5.2. Let  $I \subsetneqq J \text{ be ideals in a ring } R $. Then \\
	         
	         (i) If $ I^n = J^n$  for some $ n \geqq 2$ , $ I J^{n-1} = J^n$, i.e. $I$ is a reduction of $J$. Thus, in particular, if an ideal $J$ in a ring $R$ admits no proper reduction, $J$ is a big ideal. Further, $ I J^{n-1} = J^n \text { for } n \geqq 2$ does not imply $ I^n = J^n$. \\  
	         
	         (ii) If $ I^n = J^n$  for some $ n \geqq 1$, then $ I^m = J^m \text{ for all } m \geqq n$.  \\  
	         
	         (iii) If $R$ is Noetherian, then $ I^n = J^n$  for some $ n \geqq 1$ if and only if $ R[Jt]/R[It] $ is a finitely generated $R-$module.\\            
	         
	         Proof. (i) If $I^n = J^n$  for some $ n \geqq 2$,  $$ IJ^{n-1} \subset J^n = I^n$$ $$ \Longrightarrow IJ^{n-1} \subset I I^{n-1} \subset I J^{n-1}.$$ Consequently $ I J^{n-1} = I^n = J^n$. Thus if an ideal $J$ admits no proper reduction, it is a big ideal. To see the last part of the statement, let $R = K[X,Y]$ be the polynomial ring in two variables $X,Y$ over a field $K$. Consider $ J = (X^3, XY, Y^4) \text{ and } I= (XY, X^3 + Y^4)$. Then it is easy to see that $ I J = J^2, \text{ but }  I^2 \neq J^2 $. \\
	         
	         (ii) In case $ n = 1$, the result is clear. However, if $ n \geqq 2$, then from (i) we have $ I J^{n-1} = J^n$. We shall, now, prove the result by induction. Note that $ I^{n+1} = I J^n= I J^{n-1} J = J^n J = J^{n+1}$. Now, by induction it is immediate that $I^m = J^ m \text{ for all } m \geqq n$. \\
	         
	         (iii) The proof the statement is clear using (ii).\\

	          If a ring satisfies the property (P), it has no non-trivial nilpotents. Thus all rings will be assumed reduced. Further, the property (P) will hold whenever $I\subsetneqq J$ implies   $I^2\subsetneqq J^2$ since if $I\subsetneqq J$ is a counter example then for some $ n > 1,I^n = J^n$. Choose  $n$ least such that $I^n = J^n$ then clearly  $I^{n-1}\subsetneqq J^{n-1}$ gives a counter example. Let us also note that $R$ will satisfy (P) if the property holds for all pairs of ideals in $R$ of the form $I\subsetneqq id(I,f) = J$ where $f \in R - I$. Hence $R$ satisfies the property (P) if and only if for every ideal $I$ in $R$ and $ f\notin I,$ if $ f^2 \in I^2 $ and $ fI \subset I^2$, then $f \in I$.  
	         We, now,  note that the the property does not hold even in Noetherian domains.\\        
	         
	          {\bfseries Example} 6.2. Let $ R= \mathbb{Z}[X^3, X^4, X^5], I = id.(X^3,X^4, 2X^5), J=id(X^3,X^4,X^5).$ Then $I\subsetneqq J$, but $ I^2 = J^2= id(X^6,X^7,X^8)$.\\

                {\bfseries  Lemma} 7.2. Let $I\subsetneqq J$ be ideals in an integral domain $R$ such that $I$ is invertible, then $ I^n \subsetneqq J^n$ for all $ n \geqslant 1.$\\
                
               {\bfseries  Proof}. Assume $ I^n = J^n$, then $$ I^n = I^{n-1} I \subset I^{n-1} J \subset J^{n-1} J =J^n =I^n $$  Hence  $$ I^{n-1}J = I^n$$ Thus as $I$ is invertible $ I =J$, which is not true. Consequently $ I^n \subsetneqq J^n$ for all $ n \geqslant 1$.\\
                
              {\bfseries  Corollary} 8.2. If $R$ is a Dedekind domain, then for any two non zero ideals $ I \subsetneqq J$ in $R, I^n \subsetneqq J^n$ for all $ n\geqslant 1$.\\
                
              {\bfseries  Lemma} 9.2.  Let $I\subsetneqq J$ be ideals in a ring $R$ such that $J \nsubseteqq \sqrt{I} $, then $ I^n \subsetneqq J^n$ for all $ n \geqq 1.$\\
                
               {\bfseries  Proof}. Let $ f \in J -\sqrt{I} $, then $$ I \subsetneqq id(I,f)= I_1 \subset J.$$ If $I^n = I_1^n$, then $ f^n \in I^n \subset I$. Thus $ f \in \sqrt{I}$, which is not true. Hence $ I^n \subsetneqq I_1^n \subset J^n.$ Thus the result follows.\\
                
             { \bfseries  Remark} 10.2(i) If $ \wp $ is a prime ideal in $R$  and $ \wp \subsetneqq J$  for an ideal $J$ in $R$, then $ \wp^n \subsetneqq J^n$ for all $ n \geqq 2.$\\
                
                 (ii)  From the above lemma, we note that if  $I \subsetneqq J$ be ideals in a ring $R$ then we can have $ I^n = J^n$ for some  $n > 1$ only when $ \sqrt{I} =\sqrt{J}$. However, even when  $ \sqrt{I} =\sqrt{J}$, we can have $ I^n \subsetneqq J^n$ for all $ n> 1$. e.g. if $R$ is a Dedekind domain, then for any non zero ideal $J$ in $R$ if $ I = J^n, n\geqq 2$, then  $ \sqrt{I} =\sqrt{J}$, but $ I^n \subsetneqq J^n$ for all $ n \geqq 1$. We give below another example:\\
                
               {\bfseries   Example} 11.2. Let $R = \mathbb{Z}[X , Y]/id(X^2 +Y^2)$, and  $ I = id(\bar{X}), J = id( \bar{X}, \bar{Y})$ be ideals in $R$. Clearly $ I \subsetneqq J$ and $ \sqrt{I} = \sqrt{J} = J.$ Note that $ I^n \subsetneqq J^n$ as $ \bar{X}^{n-1} \bar{Y}\notin id(\bar{X}^n), $ since otherwise $ X^n - X^{n-1}Y = ( X^2 +Y^2)h(X,Y)$ for some $h(X,Y) \in \mathbb{Z}[X , Y]$. This,however, is not true as putting $X = iY$ in this equation gives $ X^n + iX^n = 0$.\\
                
                 {\bfseries Theorem} 12.2. Let $R$ be an integral domain which is not a field. Then there always exist ideals $ I \subsetneqq J$ in $R[X]$ such that $ I^2 = J^2$. Thus $ R[X]$ does not satisfy the property (P).\\
                 
                {\bfseries  Proof}. Let $ a\in R$ be a non-zero,non-unit element. Let $ I =id(X^4, a^3 X, a X^3, a^4)$ and $ J = I + id(a^2X^2)$ be ideals in $R[X]$. Then it is easy to check that $ I^2 = J^2$. Thus to prove our claim, it suffices to verify that $ I \subsetneqq J$ . This is true if and only if $a^2 X^2 \notin  I = id(X^4, a^3 X, a X^3, a^4).$ Assume the contrary, then $$ a^2X^2 = X^3 g_1(X) + a^3 X g_2(X) + aX^3 g_3(X) + a^4 g_4(X)$$ for some $g_i(X) \in R[X]$ for all $ i \geq 1.$ Note that, clearly this is not true as $a^2 X^2$ does not occur in any  term on the right  hand side of the equation. Hence the result holds.\\ 
                
                {\bfseries Observation } 13.2.              
                 If $I$ is a regular ideal in a Noetherian ring $R$, and $ I \neq I^* = \cup \{ (I^{l+1} : I^l) l \geq 1\}$, then $R$ does not satisfy the property (P) since $ I^{*k} = I^k$ for all large k [7, Theorem 2.1]\\  
                
                 {\bfseries Theorem} 14.2. Let $R$ be a Noetherian integral domain. Then\\
                 
                 (i) The ring $R$ satisfies the property (P) if and only if $ I^* = I $ for all ideals $I$ in $R$.\\
                 
                 (ii) If every ideal in $R$ is integrally closed then $R$ satisfies the property (P).\\
                
                {\bfseries Proof}. (i) Let $R$ satisfies the property (P). Then $ I^* = I$ for all ideals $I$ in $R$ by [7, Theorem 2.1]. Conversly, let $ I^* = I$ for all ideals $I$ in $R$. Assume  $I\subsetneqq J$ be a pair ideals in a ring $R$ for which the property (P) fails. Then  $ I^n =J^n$ for all $ n \gg$. Hence by [7, Theorem 2.1] , $ J \subset I^* = I$. This contradicts our assumption that  $I\subsetneqq J$.  Cosequently $ I^n \subsetneqq J^n$ for all $n$. Thus $R$ satisfies (P).\\ 
                
                (ii) The assretion follows by (i) and [7, Remark 2.3.3].\\
                
                In view of above, the Noetherian domains which satisfy the property (P) are precisely those in which every non-zero ideal is a Ratliff-Rush ideal. it is natural to ask if a Ratliff- Rush ideal a big ideal in a Noetherian integral domain. We give an example to show that this is not true.\\
                
                {\bfseries  Example} 15.2. Let $ R = K[[t^3,t^4]]$, where $K$ is a field and $t$ is an indeterminate over $K$. As noted in $[3,\mbox{ Example }  1.11]$ all  the powers of  maxinal ideal $ \mathbf{m}=  id(t^3,t^4)$ in $R$ are Ratliff -Rush ideals. However as $ id(t^6,t^7) = t^3 \mathbf{m} \subsetneqq \mathbf{m}^2 = id(t^6,t^7,t^8)$, but $ (t^3 \mathbf{m})^2 = \mathbf{m}^4 = id( t^{12},t^{13},t^{14})$, the ideal $ \mathbf{m}^2$ is not a big ideal.\\

                 As proved in $[7, Prop. 3.1]$, we also have.\\
                
                {\bfseries Theorem} 16.2. Let $R$ be a Noetherian integral domain. If $R$ satisfies the property (P), then dimenstion of $R$ is $ \leqq 1$.\\
                
                {\bfseries  Proof}. Assume the dimenstion of $R$ is $ \geqq 2$. Then there exist $ x,y\in R$ such that $ id(x,y)$ has  height 2. By $[5, Chap. V, Thm. 4.14] $, $ \{ x, y\} $ is an independent set. Put $ I = id (x^4, x^3 y, x y^3, y^4))$. As $ x,y$ are independent $ x^2 y^2 \notin I $. Thus $  I \subsetneqq J = id (I, x^2 y^2)$, but $ I^2  = J^2$. Hence the assertion follows.\\
                
                {\bfseries   Theorem} 17.2 Let $ I \subsetneqq J $ be ideals in a Noetherian ring $R$. Then $I^n \subsetneqq J^n$ for all $ n \geqq 1$ if either of the conditions (i) $\text{height} (I)\neq \text{ height}  (J)$ or (ii) $\text{grade} (I) \neq \text{ grade}(J)$ or (iii) set of minimal primes of $I  \neq $ set of minimal primes of $J$ or (iv) $\text{dim} (R/I) \neq \text{dim}(R/J)$ or (v) $ \text{radical}(I) \neq \text{radical}(J)$ hold.\\
                
                 {\bfseries  Proof}. (i) This follows since for any ideal $ K$ in the ring $R,\text{height}(K) = \text{height} (K^n) $ for all $ n\geqq 1$.\\
                 
                 (ii) By [8, lemma 3.2], the grade of an ideal is equal to  the grade of its radical. Hence the assertion follows. \\
                 
                 (iii) Clearly, the set of minimal primes of an ideal is equal to the set of minimal primes of any its powers. Thus the result follows.\\
                 
                 (iv) For any ideal $K$ in $R,\text{dim} (R/K) = \text{dim} (R/K^n)$ for any $ n \geqq 1$. Hence the result follows. \\
                 
                 (v) As radical of an ideal is equal to the radical of any its power, the result follows.

                 \section{Big Ideals}
              
               In this section, we shall study the existance of big ideals in rings. Let us note that :\\
               
               O.1.3. Big ideals are invariant under isomorphism.\\
               
               O 2.3. Let a ring $S$ be a faithfully flat extension of a ring $R$. If for an ideal $J \subset R, JS$ is a big ideal of $S$ , then $J$ is a big ideal in $R$.\\            
             
             {\bfseries Theorem} 3.3. Let $R$ be a Noetherian integral domain then an ideal $J$ in $R$ is a big ideal if and only if whenever $ I \subsetneqq J$ is an ideal, then $ J \nsubseteqq I^*$.\\
  
              {\bfseries Proof}. Let $J$ be an ideal in $R$ such that whenever $ I \subsetneqq J$ is an ideal, then $ J \nsubseteqq I^*$. If $J$ is not a big ideal then there exists an ideal  $I_1 \subsetneqq J$ such that $ I_1^n = J^n $ for some $ n > 1$. Then for any $ r \leq n$, we have $$ I_1^r J^{n-r} \subset J^n = I_1^n = I_1^r I_1^{n-r} \subset I_1^r J^{n-r}$$ Hence $ I_1^r J^{n-r} = I_1^n = J^n$ for all $ r \leq n$. Thereforec $$ I_1^{n+1} \subset J^{n+1} =J I_1^n = J I_1^{n-1} I_1 = I_1^n I_1 = I_1^{n+1}$$ Hence $ I_1^{n+1} = J^{n+1}$. Consquently, by induction $ I_1^m = J^m $ for all $ m \geq n$. Thus $ J^* = I_1^* $ and hence $ J \subseteqq I_1^*$. This cotradicts our assumption that   $J \nsubseteqq I_1^*$. Hence $J$ is a big ideal. Conversely let $J$ be a big ideal, and let $ I_1 \subsetneqq J$ be an ideal. If $ J \subseteqq I_1^*$, then by [7, Theorem 2.1], $ J^t \subseteqq  I_1^{*t}= I_1^t \subseteqq J^t$ for all $ t $ large. Hence $ J^t = I_1^t$ for all $t$ large. This cotradicts the assumption that $J$ is a big ideal. Hence the result folows.\\
              
               {\bfseries Remark} 4.3.  In a Noetherian integral domain $R$, an ideal $J$ is a big ideal if and only if whenever $ I \subsetneqq J$ is an ideal, then  $ J I^n \nsubseteqq I^{n+1}$ for all $ n\geq 1.$\\
                   
               { \bfseries Lemma} 5.3. Let $J$ be a big ideal in an integral domain $R$. Then for any invertible ideal $ \textit{A} \subsetneqq R, \textit{A}J$ is a big ideal.\\ 
                
                  Pf. Let $I \subsetneqq \textit{A}J $ be an ideal in $R$. Then
                
                                     $$  \textit{A}^{-1} I \subsetneqq J$$
                                $$ \Longrightarrow    ( \textit{A}^{-1} I)^n \subsetneqq J^n$$  $$ \Longrightarrow I^n \subsetneqq (\textit{A} J)^n .$$ Hence the assertion follows.\\  
                                                                	
            {\bfseries  Remark} 6.3.  Any invertible ideal in $R$ is a big ideal. This follows by taking $ J=R$ in the Lemma ( As an exception we are considering $R$ an ideal).\\

               \textbf{Lemma} 7.3.Let $R$ be a Noetherian integral domain. An ideal $J \subset R $ is a big ideal if and only if there  does not exist any ideal $ I  \subsetneqq J$ such that $ I^* = J^*$.\\
               
             \textbf{ Proof}. Assume $J$ is a big ideal. If for an ideal $ I  \subsetneqq J, I^* = J^*$, then since by [7, Theorem 2.1], $I^{*n} = I^n \mbox{ and } J^{*n} = J^n $ for all n large, we get $ I^n = J^n$ for all n large. This contradicts that $J$ is a big ideal. The converse is also clear by [7, Theorem 2.1 ]   \\

              {\bfseries  Theorem} 8.3. Any maximal ideal in a regular  ring is a big ideal.\\
               
               Proof. Let $ \textbf{M}$ be a maximal ideal in a  regular ring $R$, and let $\textbf{I}\subsetneqq \textbf{M} $ be an ideal. If $ \textbf{M} \neq \sqrt{I}$, then $ \textbf{I}^n \subsetneqq \textbf{M}^n$ for all $ n \geqq 1$ by Lemma 8
               .2. Now, let $ \textbf{M} = \sqrt{I} $. Then $ \textbf{I} $ is $ \textbf{M}$ primary ideal. Hence $\textbf{I}_{\textbf{M}} \subsetneqq \textbf{M}R_{\textbf{M}}$ is  $ \textbf{M}R_{\textbf{M}}-$ primary ideal in $R_{M} $. As $R_{M} $ is regular local ring $ \textbf{M}R_{\textbf{M}}$ is generated by a system of paramers. Hence by [4, Corollary 2.4] $\textbf{M}R_{\textbf{M}}$ is basic. Thus $(\textbf{I}_{\textbf{M}})^n \subsetneqq    (\textbf{M}R_{\textbf{M}})^n$ since othewise $\textbf{I}_{\textbf{M}}$ will be a proper reduction of $\textbf{M}R_{\textbf{M}}$, which contradicts that $\textbf{M}R_{\textbf{M}}$ is basic. 
               Consequently $ \textbf{I}^n \subsetneqq \textbf{M}^n$ for all $ n \geqq 1$ i.e., $ \textbf{M}$ is a big ideal. \\

              {\bfseries  Remark} 9.3.  (i) Let $\textbf{Q}$ be an $\textbf{M}-$ primary ideal  in a Noetherian local domain  $ (R , \textbf{M})$ . Then if $ \textbf{Q} $ is generated by a system of parameters, it  is a basic ideal  by [4, Corollary 2.4].   Thus if $ \textbf{I} \subsetneqq \textbf{Q}$ , then  $ (\textbf{I})^n \subsetneqq ( \textbf{Q})^n$ for all $ n \geqq 1$ since otherwise $ \textit{I}$ wiil be a proper reduction of $ \textbf{Q}$. \\   
                          
              (ii)  Let $ M_1. \ldots, M_n$ be distinct maximal ideals in a regular ring $R$  Then by Theorem 9.3 and [4, Theorem 3.6],  $ I = M_!.M_2.\ldots M_n$ is big ideal.\\
              
              {\bfseries Theorem }10.3 (i) Let $J$ be an ideal in a Noetherian ring $R$ generated by an $R-$sequences. Then $J$ is a big ideal.\\
              
              (ii) Let $K$ be a field, then every prime ideal in $ K[X, Y]$ is a big ideal. Moreover, every ideal generated by two elements in $ K[X_1, \ldots, X_n], n\geqq 2, $ is a big ideal.\\
                 
                 Proof. (i) Let $J = id.(a_1,\cdots,a_n)$ where $ \{a_1, \cdots, a_n\}$ is an $R-$sequence. By [9, Theorem 2.1], $J$ is an ideal of principal class, Hence by [2, Corollary 2.4], $J$ is basic. Consequently $J$ is a big ideal by Lemma, 8.3.\\
                 
                 (ii By Theorem 8.3 and  Remark 6.3, it is immediate that every prime ideal in $K[X, Y]$ is a big ideal. Further, let $ I= id(f,g) $ be any ideal in $ K[X_1, \ldots, X_n], n\geqq 2$. If $ I $ is principal, it is clearly big. Now assume, $I$ is not principal. Thus $ I= h .id(f_1,g_1) $ where $(f_1, g_1) = 1$  i.e., g.c.d. of $f_1 \text{  and } g_1 $ is a unit. Then $ \{f_1, g_1\}$ is a regular sequence. Thus $ I $ is a big ideal by Theorem 9.3 and Remark 6.3. \\

                   {\bfseries Remark} 11.3 Let  $R[ X_1, \cdots, X_n]$ be a polynomial ring over a Noetherian ring $R$. Then for $ k \leqq n $  the ideal  $ id( X_1-a_1, \cdots, X_k -a_k), \mbox{ where } a_i \in R \mbox{ for all } i=1,2 \cdots , k$, is a big ideal.\\
                 
                   {\bfseries Theorem} 12.3.(i) Let $I$ be an ideal in a Noetherian ring $R$. Then $I$ is a big ideal if $I_M$ is a big ideal in $R_M$ for every maximal ideal $M \supseteq I$ in the ring $R$.\\
                   
                   (ii) Let $R$ be an almost Dedekind domain. Then every ideal in $R$ is a big ideal.
                   
                   Proof.  (i) Assume $I_M$ is a big ideal in $R_M$ for every maximal ideal $M \supseteq I$ in the ring $R$. Let $ J \subsetneqq I$ be any ideal. Then there exists a maximal ideal $M \supsetneqq I$  such that $ J_M \subsetneqq I_M$. Hence, as $I_M$ is a big ideal, $ (J_M)^n \subsetneqq (I_M)^n $ for all $ n \geqq 1$. Consequently $ J^n  \subsetneqq I^n $ for all $ n \geqq 1$. Thus $I$ is a big ideal.\\ 
                   
                   (ii) Let $I \subsetneqq J$ be ideals in $R$. Then there exists a maximal ideal $M \supset J$ such that $ I_M \subsetneqq J_M$. Hence by the Lemma 7.2, $ I^n_M \subsetneqq J^n_M$. Consequently $I^n \subsetneqq J^n$ for all $ n \geqq 1$. \\ 
                   
                   {\bfseries Theorem} 13.3 Let $R$ be a regular ring and $ \wp $ be a prime ideal in $R$. If for an ideal $I \subsetneqq \wp, I_{\wp} \neq \wp R_{\wp}$. Then $ I^n \subsetneqq \wp^n \text { for all } 
                  n \geq 1$.\\
                  
                 {\bfseries  Proof.} In view of Lemma 8.2, it is enough to prove the assertion in case $ \sqrt{I} = \wp$. In this case, $ I_{\wp} \text { is } \wp R_{\wp}-$primary. As $R$ is regular $ R_{\wp}$ is a regular local ring. Hence $\wp R_{\wp}$ is basic, and consequently it is a big ideal. Therefore $ I^n_{\wp} \subsetneqq \wp^nR_{\wp}$, and hence $ I^n \subsetneqq \wp^n$ for all $ n \geqq 1$. \\   
                  
                  { \bfseries Remark.} 14.3 (i)  Note that in the Theorem $ I_{\wp} =\wp R_{\wp} $ if and only if there exists $ s\notin \wp$ such that $ s\wp \subset  I$. Thus $ \wp \subset Z_R(R/I)$. Hence if $ \wp \nsubseteq Z_R(R/I)$, then $ I^n \subsetneqq \wp^n \text { for all } n \geq 1$.  \\
                  
                  (ii) Let $R$ be an integral domain, then even if $I_{\wp} = \wp R_{\wp}$, we can have $ I^n \subsetneqq \wp^n \text { for all } n \geq 1$ e.g. take $ I = s\wp$  for some $s \notin \wp$. For this part we need only to assume that $\wp$ contains a non-zero divisor.\\  
                  
                  We, now, prove two results which fall short of showing that an integral domain $R$ in which every finitely generated ideal is a big ideal is integrally closed.\\

                   {\bfseries Theorem} 15.3 Let $R$ be an integral domain in which every ideal generated by atmost three elements is a big ideal. Let for non- zero elements $ x,y \in R, x \notin Ry$.  Then if  $ x^2 \notin Ry, x/y $ is not integral over $R$.\\
                   
                   Proof. Let $K$ be the field of fractions of $R$. As $ x \notin  Ry, x/y \in (K-R)$.  If $ x/y$ is integral over $R$, then $$  (x/y)^n + c_1 (x/y)^{(n-1)}+ \ldots + c_n =0 $$  	 where $ n > 1$ and $ c_{i's} \in R.$  $$\Longrightarrow x^n + c_1 x^{n-1}y + \ldots + c_n y^n = 0 $$ $$ \Longrightarrow  x^n \in (x,y)^{n-1} y  $$ $$\Longrightarrow  y (x,y)^{n-1} = (x,y)^n   $$  $$  \Longrightarrow  y^t (x,y)^{n-1} = (x,y)^{t+n-1}$$ for all $ t \geqq 1$.   Hence $$ y^{n-1} (x,y)^{n-1} = (x,y)^{2(n-1)}$$ As every ideal generated by atmost three elements in $R$ is big, we conclude $$ y(x,y) = (x,y)^2 = ( x^2, y^2, xy)$$  Thus $ x^2 \in (y)$.A contradiction to assumption.  Hence $ x/y$ is not integral over $R$.\\
                   
                     {\bfseries Theorem} 16.3 Let $R$ be an integral domain with quotient field $K$ in which every finitely generated ideal is a big ideal. Then if  $ x/y \in (K-R)$, either $ x/y \mbox { or } y/x^2$ is not integral over $R$. \\
                     
                     Proof. Let $z \in R$ be a non-zero, non-unit. Then $  z, z^2 \notin Rz^3.$ Hence by Theorem 13.3, we get $ 1/z^2 = z/z^3$ is not intergral over $R$. Therefore $ 1/z$ is not integral over $R$. Now, let $ x/y \in (K- R)$. Then $ x \notin Ry$. By Theorem 12.3, if $ x^2 \notin Ry$, then $ x/y$ is not integral over $R$. Assume $ x^2 \in Ry$.  Then $ x^2 = \mu y$ for an element $ \mu \in R$. Clearly $x$ is not a unit since otherwise $y$ is also a unit, and thus $ x/y \notin (K-R)$. If $\mu$ is unit then  $ x/y = 1/\mu^{-1}x$ is not integral over $R$ as seen above. However if $\mu$ is not a unit then $ y/x^2 = 1/ \mu$ is not integral over $R$.\\ 
                     
                     {\bfseries Theorem.} 17.3 Let $R$ be an  integral domain( not necessarily Noetherian). If every finitely generated ideal in $R$ is a big ideal,then $R$ has no regular sequence of length $\geq 2$.\\
                     
                     Proof. Let $ \{a,b\}$ be a regular sequence in $R$. We shall first show that $ ab \notin (a^2, b^2 )$. If not, then $$ ab = \lambda a^2 + \mu b^2$$ for some $ \lambda, \mu \in R$.  $$ \Longrightarrow \mu b^2 \in Ra$$ $$\Longrightarrow \mu \in Ra$$ Let $ \mu = pa$. Then $$ ab = \lambda a^2 + pab^2$$ $$ \Longrightarrow b = \lambda a + p b^2$$ $$ \Longrightarrow b( 1- pb) = \lambda a$$ $$\Longrightarrow 1-pb \in Ra$$ $$ \Longrightarrow 1 \in (a,b).$$ This is not true. Hence $ ab \notin (a^2,b^2). $ Thus $ (a^2,b^2) \subsetneqq (a,b)^2$. Note that $ (a,b)^2 (a^2,b^2) = (a,b)^4$. Hence $(a^2,b^2)$ is a reduction of $ (a,b)^2$. Thus by Lemma 5.2. $(a,b)^2$ is not a big ideal. This contradicts our assumption on $R$. Hence the result follows.\\

                    \section{ Exponentially Equal Ideals} If an ideal $J$ in a ring $R$ is not a big ideal, then there exists an ideal $I \subsetneq J$ such that $ I^n = J^n$ for all $ n \gg$. Based on this fact, we define the concept of exponentially equal ideals and  prove that in a local ring $R$ if $ I \subset J $ are ideals which are exponentially equal then there exists an ideal $ I^{'} \subset I $ minimal with respect to the property that $I^{'}$ is exponentially equal to the ideal $J$. The proof follows the arguments  in [6], used for the existance of minimal reductions of an ideal. \\
                 
                { \bf Definition.}1.4 Let  $R$ be a ring, and $ I , J$ be ideals in $R$. We shall say that $I$ is  exponentially equal to $J$ if $ I^n = J^n $ for all $n \gg$.	\\  
                 
                {\bf Remark} 2.4. (i) Exponentially equal ideals have same set of minimal prime ideals, and have same multiplicity at every common minimal prime. \\
                 
                 (ii) No two ideals in a Prufer domain are exponentially equal [1, Exercise 1, page 284]\\
                 
                 (iii) Exponential equality is an equivalence relation on ideals.\\
                 
                 {\bf Lemma. }3.4 If $I$ is a regular ideal in a Noetherian ring $R$, then for an ideal $J \subset R$, $I$ is exponentially equal to $J$ if and only if $ I^* = J^*$.\\
                 
                 {\bf Proof}. The proof is immediate from [7, Theorem 2.1 ]. \\   
                            
               {\bf Theorem} 3.5  Let $ R $ be a Noetherian local ring with maximal ideal $ \mathbf{m}$, and let $ I \subsetneqq J$ be  ideals in $R$. Then  \\
                 
                 (i) For any $ n \geq 1$, $ I^n = J^n$  if and only if  $ (I + \mathbf{m}J) ^n = J^n$. \\
                 
                 (ii) If $ I^n = J^n$ for all $ n \gg$, then  there exists an ideal $ I_2 \subset I$ minimal with the property $(I_2)^n = J^n$ for all $ n \gg $.\\
                 
                {\bf Proof.}  (i)  Let $ I^n = J^n$. Then  $ J^n = I^n \subset (I + \mathbf{m} J)^n \subset ( I^n + \mathbf{m} J^n ) \subset J^n  $. Hence $(I + \mathbf{m} J)^n = J^n$. Conversely if $(I + \mathbf{m} J)^n = J^n$, then $ J^n \subset I^n + \mathbf{m} J^n \subset J^n$. Therefore $ I^n + \mathbf{m} J^n = J^n$, which implies $ I^n = J^n$. Thus (i) follows.  \\
                 
                 (ii) Let $$  \sum = \{  \mathbf{K} \subset I, \text{ an ideal }\mid  \mathbf{K}^n = J^n  \mbox { for all } n \gg \}.$$  Clearly, $ \sum \neq \phi \mbox { since } I \in \sum  $.  By (i), for any $ \mathbf{K} \in \sum, \mathbf{K} + \mathbf{m}J \in \sum$. Moreover, $$ ( \mathbf{K} + \mathbf{m}J)/ \mathbf{m}J \subset J/\mathbf{m}J$$
                 for every $\mathbf{K} \in \sum .$ As $ dim_{R/\mathbf{m}} (J/\mathbf{m} J) < \infty $, there exists $I_1 \in \sum$ such that 
                 $ dim_{R/\mathbf{m}}(I_1 + \mathbf{m}J ) / \mathbf{m} J $ is least. Choose $ x_i \in I ; i= 1,2,\cdots, t$ such that $ x_i + \mathbf{m}J ; i=1,2,\ldots, t$ generate the $ R/\mathbf{m}-$ vector space $ ( I_1 + \mathbf{m}J )/ \mathbf{m} J $. Put $I_2 = id.(x_1,\ldots, x_t)$, then $ I_2 + \mathbf{m} J = I_1 + \mathbf{m} J $ . Further, if $\sum_{1}^{t}\alpha_i x_i  \in \mathbf{m} J$ where $  \alpha_i \in R $, then $ \alpha_i \in  \mathbf{m}$. Hence $$ I_2 \cap \mathbf{m} J  \subset  \mathbf{m} I_2$$ Note that, as $ I_1 + \mathbf{m} J  = I_2 + \mathbf{m} J , ( I_1 + \mathbf{m} J)^n  = (I_2 + \mathbf{m})^n = J^n . $ Now, let $ I_3 \subset I_2$ be such that $( I_3)^n = J^n$, then $( I_3 + \mathbf{m} J)^n = J^n$. Moreover $$ (I_3 + \mathbf{m} J)/\mathbf{m}J  \subset   (I_2 + \mathbf{m} J)/\mathbf{m}J =(I_1 + \mathbf{m} J)/\mathbf{m}J  \subset    J/\mathbf{m}J$$ By choice of $ I_1$, it follows that $ I_3 + \mathbf{m} J = I_2 + \mathbf{m} J = I_1 + \mathbf{m} J$. Now, if $ \lambda \in I_2$, then $ \lambda = x + b$ where $ x \in I_3, b \in \mathbf{m} J$. Thus $ \lambda - x \in I_2 \cap \mathbf{m} J \subset \mathbf{m} I_2$. Hence $ I_2 \subset I_3 + \mathbf{m} I_2$. Consequently $ I_2 \subset I_3$. This implies $ I_2 = I_3$. Hence $I_2$ is minimal.\\    
                 
                It would be interesting to know the answer of the following:\\
                 
                {\bf  Question} 1.Let $K$ be a field. Is every prime ideal in the polynomial ring $ K[X_1, \ldots, X_n], n\geqq 3,$ a big ideal?\\

                {\bf Question} 2.  Let $K$ be a field. Can we characterize big ideals in  the polynomial ring $ K[X_1, \ldots, X_n], n\geqq 2$?\\

                 \begin{center}
                 	\bf{ACKNOWLEDGEMENT} \end{center} 
                 I am thankful to Melvin Hochster for some useful e-mail exchanges and to Google  for access to literature.\\

                             \begin{center}
                               	 \bf{REFERENCES} \end{center} 
                     1. R.Gilmer, Multiplicative Ideal Theory, Queen's Papers in Pure and Applied Math., No.12, 1968.\\
                     
                     2. James H. Hays, Reductions of Ideals in Commutative Rings, Trans. of Amer. Math. Soc., Vol.177, pp. 51-63.  \\         
                                
                     3. William Heinzer, David Lantz and Kishor Shah, The Ratliff- Rush Ideals in a Noetherian Ring, Communications in Algebra, 1992. \\
                     
                     4. James H. Hays, Reduction of Ideals in Commutative Rings, Trans. of American Math. Soc., Vol. 177, 1973, pp. 51-63.\\
                     
                     5. Ernst Kunz, Introduction to Commutative Algebra and Algebraic Geometry, Birkhauser, 1985.\\
                     
                     6. D.G. Northcott and D.Rees, Reductions of ideals in Local Rings, Proc. of The Cambridge Phil. Soc., Vol. 50, 1954.\\  
                 	
                 	 7. L.J.Ratliff and David E. Rush, Two notes on reduction of ideals, indiana University mathematics Journal, Vol. 27, No.6(1978) \\
                 	 
                 	 8. D.Rees, A Theorem of Homological Algebra, Proc. of The Cambridge Phil. Soc., Vol. 52, 1956.\\
                 	 
                 	 9. D.Rees, The Grade of an Ideal or Module, Proc. of The Cambridge Phil. Soc., Vol. 53, 1957.\\
                 	 
                 	 10. Pramod K. Sharma, Powers Vs. Powers  (In preparation)

          \end{document}